\documentclass[11pt,reqno]{amsart}
\usepackage{amssymb,amsthm,amsmath,enumerate}
\usepackage[all]{xy}

\newcommand{\myemph}{\emph}

\newcommand{\defeq}{=_{\mathrm{def}}}
\newtheorem{theorem}{Theorem}[section]
\newtheorem*{theorem*}{Theorem}
\newtheorem{lemma}[theorem]{Lemma}
\newtheorem{proposition}[theorem]{Proposition}
\newtheorem{corollary}[theorem]{Corollary}

\theoremstyle{definition}
\newtheorem{definition}[theorem]{Definition}

\theoremstyle{remark}
\newtheorem{remark}[theorem]{Remark}

\newtheorem{example}[theorem]{Example}

\newcommand{\Implies}{\Rightarrow}

\newcommand{\Iff}{\Leftrightarrow}
\renewcommand{\iff}{\Leftrightarrow}

\newcommand{\lscott}{ \{ }
\newcommand{\rscott}{ \} }

\newcommand{\catE}{\mathcal{E}}

\newcommand{\mapS}{\mathcal{S}}
\newcommand{\yon}{\mathbf{y}}
\renewcommand{\hom}{\operatorname{Hom}}

\newcommand{\im}{\operatorname{Im}}
\newcommand{\myker}{\operatorname{Ker}}

\newcommand{\iso}{\cong}

\newcommand{\cov}{\operatorname{Cov}}

\newcommand{\Psh}[1][]{\operatorname{Psh}_{#1}}
\newcommand{\Sh}[1][]{\operatorname{Sh}_{#1}}
\newcommand{\JSh}[1][]{\catE_{J}}
\newcommand{\smallpower}{\mathcal{P}}
\newcommand{\sub}[1][]{\operatorname{Sub}_{#1}}
\newcommand{\ssub}[1][]{\operatorname{Sub}_{\mapS}}

\newcommand{\myomega}{\Omega_\mapS}

\newcommand{\cl}{\mathbf{C}}
\newcommand{\inbar}[1][]{\, \overline{\in}_{#1} \,}
\newcommand{\nibar}[1][]{\, \overline{\ni}_{#1} \,}

\newcommand{\peejay}[1][]{\mathcal{P}_{J_{#1}}}
\newcommand{\inj}{\, \in_J \, }
\newcommand{\nij}{\, \ni_J \, }
\newcommand{\ssubj}{\operatorname{Sub}_{\mapS_J}}
\newcommand{\asf}{\mathbf{a}}
\newcommand{\incl}{\mathbf{i}}

\SelectTips{cm}{}
\newdir{ >}{{}*!/-7pt/@{>}}

\begin{document}

\author{S. Awodey}
\address{Department of Philosophy, Carnegie Mellon University}
\email{awodey@cmu.edu}

\author{N. Gambino} 
\address{Department of Computer Science, University of Leicester}
\email{nicola.gambino@gmail.com}

\author{P. L. Lumsdaine} 
\address{Department of Mathematical Sciences, Carnegie Mellon University}
\email{plumsdai@andrew.cmu.edu}

\author{M. A. Warren}
\address{Department of Mathematics and Statistics, University of Ottawa}
\email{mwarren@uottawa.ca}

\title[Lawvere-Tierney sheaves in Algebraic Set Theory]{Lawvere-Tierney sheaves \\ in Algebraic Set
Theory} 

\date{October 22nd, 2008} 

\begin{abstract}
We present a solution to the problem of defining a counterpart
in Algebraic Set Theory of the construction of internal
sheaves in Topos Theory. Our approach is general
in that we consider sheaves as determined by 
Lawvere-Tierney coverages, rather than by Grothendieck
coverages, and assume only a weakening of the axioms for 
small maps originally introduced by Joyal and Moerdijk,
thus subsuming the existing topos-theoretic results. 
\end{abstract}

\maketitle

\section*{Introduction}

Algebraic Set Theory provides a general framework for the study of 
ca\-te\-go\-ry-theoretic models of set theories~\cite{JoyalA:algst}. The 
fundamental objects of interest are 
pairs $(\catE, \mapS)$ consisting of a category~$\catE$ 
equipped with a distinguished family of maps $\mapS$, whose
elements are referred to as small maps. The category~$\catE$ 
is thought of as a category of classes, and $\mapS$ as 
the family of functions between classes whose fibers are sets. 
The research in the area has been following two general directions:
the first is concerned with isolating axioms for the pair $(\catE,
\mapS)$ that guarantee the existence in $\catE$ of a model for a
given set theory; the second is concerned with the
study of constructions, such as that of internal sheaves, 
that allow us to obtain new pairs $(\catE, \mapS)$ from given ones, 
in analogy with
the existing development of Topos Theory~\cite[Chapter~5]{MacLaneS:shegl}.
The combination of these developments is intented to give general 
methods that subsume the known techniques to define sheaf and realizability
models for classical, intuitionistic, and constructive set theories~\cite{FourmanM:shems,FreydP:axic,GambinoN:heyics,GraysonR:forisw,LubarskyR:indrac,McCartyD:rearm,RathjenM:reaczf}.

Our aim here is to contribute to the study of the construction of internal
sheaves in Algebraic Set Theory. The starting point of our development is 
the notion of a \myemph{Lawvere-Tierney coverage}. If our ambient category 
$\catE$ were an elementary topos, Lawvere-Tierney coverages would be in 
bijective correspondence with Lawvere-Tierney local operators on the 
subobject classifier of the topos. However, since $\catE$ is assumed 
here to be only a Heyting pretopos, we work with the more general 
Lawvere-Tierney coverages. As we will see, when $\catE$ is a category 
of internal presheaves,
these correspond bijectively to the Grothendieck coverages considered
in~\cite{GambinoN:asssft}. Therefore, our development gives as a 
special case a treatment of the construction of internal sheaves 
relative to those Grothendieck sites, extending the 
results in~\cite{GambinoN:asssft}. 

Given a Lawvere-Tierney coverage, we can define an associated universal closure operator on the subobjects 
of $\catE$. While in the topos-theoretic context this step is essentially straightforward,  within our
setting  it  requires an application of the Collection Axiom for small maps. Once a closure operation is
defined, we can define the notion of internal sheaf as in the standard topos-theoretic context. Our main result asserts that the category of internal sheaves for a Lawvere-Tierney coverage is a Heyting pretopos and that it can be  equipped with a family of small maps satisfying the same axioms that  we assumed on the small maps $\mapS$ in the ambient category $\catE$. 
The first part of this result involves the definition of an associated sheaf functor, a finite-limit preserving left adjoint to the inclusion of sheaves into $\catE$. Even if our proof follows the spirit of the topos-theoretic argument due
to Lawvere~\cite[\S V.3]{MacLaneS:shegl}, there are two main differences. First, since the argument involves
the construction of power-objects, which in our setting classify indexed families of small subobjects, our proof requires 
a preliminary analysis of \myemph{locally small maps}, which form the family of  small maps between sheaves. We will apply this analysis also to prove the second part of our main result, which involves the verification that locally small maps between sheaves satisfy the axioms for a family of small maps. Secondly, we need to avoid using exponentials, since we do not assume them as part of the structure for the ambient category~$\catE$.

In recent years, there has been substantial work  devoted to isolating axioms 
on~$(\catE, \mapS)$ that provide a basic setting for both directions 
of research mentioned above. Let us briefly consider two such possible  
settings. 
The first, to which we shall refer as the \myemph{exact setting}, involves 
assuming that~$\catE$ is a Heyting pretopos, and that~$\mapS$ satisfies a 
a weakening of the axioms for small maps introduced 
in~\cite{JoyalA:algst}. The second, to which we shall refer as 
the~\myemph{bounded exact setting}, 
involves assuming that $\catE$ is a Heyting category, that $\mapS$ 
satisfies not only the axioms for small maps of the exact setting, 
but also the axiom asserting that for every object  $X \in \catE$, the 
diagonal 
$\Delta_X : X \rightarrow X \times X$ is a small map, that universal
quantification along small maps preserves smallness of monomorphisms,
and finally 
that $\catE$ has quotients of bounded equivalence relations, that is to 
say equivalence relations given by small 
monomorphisms~\cite{vandenBergB:exac}. 
Categories of ideals provide
examples of the bounded exact setting~\cite{AwodeyS:relttst}.
The exact completion and the bounded exact completion of syntactic
categories of classes arising from constructive set theories provide 
other examples of the exact setting and of the bounded exact setting,
respectively~\cite[Proposition 2.10]{vandenBergB:exac}. 
Neither setting is included in the other, and they are somehow incompatible.
Indeed, if we wish to avoid the assumption that every equivalence relation 
is given by a small monomorphism,
which is necessary to include constructive set theories~\cite{AczelP:notcst} 
within the general
development, it is not possible to assume both that $\catE$ is exact and
that every object has a small diagonal. Each setting has specific 
advantages.
On the one hand, the assumption of exactness of $\catE$ is useful to define 
an internal version of the associated sheaf functor~\cite{GambinoN:asssft}.
On the other hand, the assumption that diagonals are small has been applied 
in the coalgebra construction for cartesian 
comonads~\cite{WarrenM:catc} and to establish
results on $W$-types~\cite[Proposition~6.16]{vandenBergB:exac}. 

The choice of developing our theory within the exact setting is motivated 
by the desire for the theory to be appropriately general. Even for
the special case of Grothendieck sites, the assumption that the 
ambient category is exact seems to be essential in order to define the 
associated sheaf functor without additional assumptions on the
site~\cite{GambinoN:asssft}. In the bounded exact setting, Benno van 
den Berg and Ieke Moerdijk have recently announced a result concerning internal
sheaves on a site~\cite[Theorem 6.1]{vandenBergB:uniaas}, building on 
previous work of Ieke Moerdijk and Erik Palmgren~\cite{MoerdijkI:typttc}
and Benno van den Berg~\cite{vandenBergB:shept}. Apart from
the axioms for small maps that are part of the bounded exact setting, 
this result assumes a further axiom for small maps, the Exponentiation 
Axiom, and the additional hypothesis that the Grothendieck site has a basis. 
We prefer to avoid these assumptions since, for reasons of proof-theoretic
strength,  the Grothendieck site 
that provides a category-theoretic version of the double-negation
translation cannot be shown to have a basis within categories of classes arising 
from constructive set theories~\cite{GambinoN:heyics,GraysonR:forisw}.

\section{Algebraic Set Theory in Heyting pretoposes} 
\label{sec:asthp}

 We begin by stating precisely the axioms
for small maps that we are going to work with. 
As in~\cite{JoyalA:algst}, we assume that $\catE$ is 
a Heyting pretopos, and that ~$\mapS$ is a family of maps in~$\catE$ 
satisfying the axioms (A1)-(A7) stated below.

\begin{enumerate}[({A}1)]
\item The family $\mathcal{S}$ contains isomorphisms and is
closed under composition.
\item For every pullback square of the form
\begin{equation}
\label{equ:pb}
\xymatrix{
Y \ar[r]^{k} \ar[d]_{g} & X \ar[d]^{f} \\
B \ar[r]_{h} & A }
\end{equation} 
if $f : X \rightarrow A$ is in $\mathcal{S}$, then so is
$g : Y \rightarrow B$.
\item  For every pullback square as (\ref{equ:pb}), 
if $h : B \rightarrow A$ is an epimorphism and $g : Y \rightarrow 
B$ is in $\mathcal{S}$, then $f : X \rightarrow A$ is in $\mathcal{S}$.
\item The maps $0 \rightarrow 1$ and $1 + 1 \rightarrow 1$ 
are in $\mathcal{S}$.
\item If $f : X \rightarrow A$ and $g : Y \rightarrow B$ 
are in $\mathcal{S}$, then $f + g : X + Y \rightarrow A + B$ 
is in $\mathcal{S}$.   
\item For every commutative triangle of the form
\[
\xymatrix{
  X \ar@{->>}[rr]^{h} \ar[dr]_{f} &    & Y \ar[dl]^{g} \\
                     & A  &  }
\]
where $h : X \twoheadrightarrow Y$ is an epimorphism, if 
$f : X \rightarrow A$ is in $\mathcal{S}$,  then 
$g : Y \rightarrow A$ is in $\mathcal{S}$. 
\item For map $f : X \rightarrow A$ in $\mapS$
and every epimorphism $p : P \twoheadrightarrow X$, there exists a 
quasi-pullback diagram of the form
\[
\xymatrix{
Y \ar[r] \ar[d]_{g} & P \ar@{->>}[r]^{p} & X \ar[d]^{f} \\
B \ar@{->>}[rr]_{h} & & A }
\]
where $g : Y \rightarrow B$ is in $\mathcal{S}$ and
$h : B \twoheadrightarrow A$ is an epimorphism. 
\end{enumerate} 

We refer to~(A2) as the Pullback 
Stability axiom, to~(A3) as the Descent axiom, to~(A6) as the Quotients 
axiom, to~(A7) as the Collection axiom. Our basic axiomatisation of small maps involves  one more axiom. 
In order to state it, we need some terminology. Given a 
family of maps $\mapS$ satisfying (A1)-(A7), an 
$\mapS$-\myemph{object} is an object~$X$ such that the 
unique map $X \rightarrow 1$ is in $\mapS$. For a fixed object $A \in \catE$, 
an $A$-\myemph{indexed family of $\mapS$-subobjects} is a subobject 
$S \rightarrowtail A \times X$ such that its composite 
with the first projection $A \times X 
\rightarrow A$ is in $\mapS$. We abbreviate this by  
saying that the diagram $S \rightarrowtail A \times X \rightarrow A$ is an 
indexed family of~$\mapS$-subobjects. 
Recall that, writing $\Gamma(f) : X \rightarrowtail A \times X 
\rightarrow A$ for the evident indexed family of subobjects consisting of 
the graph of a map $f : X \rightarrow A$,  it holds that~$f$ is in
$\mapS$ if and 
only if $\Gamma(f)$ is an indexed family of $\mapS$-subobjects.
Axiom (P1), stated below, expresses that 
indexed families of $\mapS$-subobjects can be classified.

\begin{enumerate}
\item[(P1)] For each object $X$ of $\catE$ there exists an 
object $\smallpower(X)$ of $\catE$, called the \myemph{power object} of $X$, 
and an indexed family of $\mapS$-subobjects of $X$
\[
\ni_X \rightarrowtail \smallpower(X) \times X \rightarrow 
\smallpower(X) \, , 
\] 
called the \myemph{membership relation} on~$X$, 
such that for any indexed family of $\mapS$-subobjects 
$S \rightarrowtail A \times X \rightarrow A$ of $X$, there exists a 
unique map $\chi_S : A \rightarrow \smallpower(X)$ fitting in a 
double pullback diagram of the form
\[
\xymatrix{
S \ar[r] \ar@{ >->}[d] & \ni_X \ar@{ >->}[d] \\
A \times X \ar[r] \ar[d] & \smallpower(X) \times X \ar[d] \\
A \ar[r]_{\chi_S} & \smallpower(X) \, . }
\]
\end{enumerate}

Writing $\ssub(X)(A)$ for the lattice of $A$-indexed families of 
$\mapS$-subobjects of $X$, axiom (P1) can be expressed equivalently
by saying that for every map $\chi : A \rightarrow \smallpower(X)$, the functions 
\[
 \hom(A, \smallpower(X)) \rightarrow \ssub(X)(A) \, , 
\]
defined by pulling back $\ni_X \rightarrowtail \smallpower(X) \times X 
\rightarrow \smallpower(X)$, are a family of bijections, natural in~$A$. In the 
following, we often omit the subscript in the membership relation. 

Our basic axiomatisation of small maps involves only 
axioms (A1)-(A7) and (P1). Therefore, when we speak of a 
\myemph{family
of small maps} without further specification, we mean a family
$\mapS$ satisfying axioms (A1)-(A7) and (P1). In this case, elements of 
$\mapS$ will be referred to as small maps, and we speak of small 
objects and indexed families of small subobjects rather than 
$\mapS$-objects and indexed families of $\mapS$-subobjects, respectively.

One can also consider an alternative axiomatisation of small maps by 
requiring, in place of (P1), the axioms of Exponentiability (S1) 
and Weak Representability (S2), stated below.
\begin{enumerate}
\item[(S1)] If $f : X \rightarrow A$ is in $\mathcal{S}$, then 
the pullback functor $f^* : \catE/A \rightarrow \catE / X$ has a 
right adjoint, which we write $\Pi_f :  \catE / X \rightarrow 
\catE/ A$.
\item[(S2)] There exists a map $u : E \rightarrow U$ in 
$\mathcal{S}$ such
that every map $f : X \rightarrow A$ in $\mathcal{S}$ 
fits in a diagram of
form
\begin{equation} 
\label{equ:rep}
\xymatrix{
X \ar[d]_f & Y \ar[l] \ar[r] \ar[d] & E \ar[d]^u \\ 
A        & B \ar@{->>}[l]^h \ar[r] & U  }
\end{equation}
where $h : B \twoheadrightarrow A$ is an epimorphism,
the square on the left-hand side is a quasi-pullback
and the square on the right-hand side is a pullback.
\end{enumerate}

The axiomatisation of small maps with (A1)-(A7) and (S1)-(S2) is a 
slight variant of the one introduced in~\cite{JoyalA:algst}. 
The only difference concerns the formulation of the Weak Representability 
axiom, which was first introduced in~\cite{vandenBergB:uniaas}.
This is a weakening of the Representability axiom 
in~\cite[Definition~1.1]{JoyalA:algst}. The weakening involves having a 
quasi-pullback rather than a genuine pullback in the left-hand
side square of the diagram in~(\ref{equ:rep}). Example~\ref{exa:classes}
and Example~\ref{exa:setoid} illustrate how the weaker form of 
representability in (S2) is the most appropriate to consider when working 
within exact categories without assuming additional axioms for small 
maps. See~\cite{AwodeyS:relttst,SimpsonA:eleacc} for other forms of 
representability. This axiomatisation is a strengthening of the one consisting 
of~(A1)-(A7) and (P1). 
On the one hand, the combination of~(S1) and (S2) implies~(P1),
since the proof in~\cite[\S I.3]{JoyalA:algst} carries over when 
Weak Representability is assumed instead of 
Representability~\cite{vandenBergB:uniaas}. On the other hand,
Example~\ref{exa:topos} shows that there are examples satisfying
(P1) but not~(S2).  Let us also recall that
(P1)~implies~(S1) by an argument similar to the usual 
construction of exponentials from power objects in a 
topos~\cite[Proposition 5.17]{AwodeyS:relttst}. \medskip

We will make
extensive use of the internal language 
of Heyting pretoposes~\cite{MakkaiM:firocl}
This is a form of many-sorted first-order 
intuitionistic logic which allows us
to manipulate objects and maps of $\catE$ syntactically.  As an illustration
of the internal language, let us recall that for any map $f : X \rightarrow 
Y$, we have a \myemph{direct image} map~$f_{!} : \smallpower(X) \rightarrow 
\smallpower(Y)$. Assuming $f : X \rightarrow Y$ to be small, there is also
an \myemph{inverse image}
map $f^* : \smallpower(Y) \rightarrow \smallpower(X)$, which
is related to~$f_{!} : \smallpower(X) \rightarrow \smallpower(Y)$ by
the internal adjointness expressed in the internal language as follows:
\begin{equation}
\label{equ:intadj}
(\forall s : \smallpower(X))
(\forall t : \smallpower(Y)) \; 
\big( f_{!}(s) \subseteq t \Leftrightarrow s \subseteq f^*(t) \big) \, . 
\end{equation}
The internal language allow us also to give a characterisation of small maps. 
Indeed, a map $f: X \rightarrow A$  is small if and only if the following 
sentence is valid:
\begin{equation}
\label{equ:logchasmall}
(\forall a : A)(\exists s: \smallpower(X))(\forall x:X) 
\big( f(x) = a \Iff x \in s  \big) \, . 
\end{equation}
The sentence in~(\ref{equ:logchasmall}) can be understood informally as 
expressing that the fibers of $f : X \rightarrow A$ are small. Formulation 
of some of the axioms for small in the internal language can be found 
in~\cite{AwodeyS:preast}. For $s : \smallpower(X)$ and a 
formula~$\phi(x)$ where $x : X$ is a free variable, we define the restricted
quantifiers by letting
\begin{eqnarray*}
(\forall x \in s) \phi(x) & \defeq & (\forall x : X) \big( x \in s \Rightarrow \phi(x) \big) \, , \\
(\exists x \in s) \phi(x) & \defeq & (\exists x : X) \big( x \in s \land 
\phi(x) \big) \, . 
\end{eqnarray*}
For an object $X$ of $\catE$,  anonymous variables of sort $X$ are denoted as~$\_ : X$. In particular, for a 
formula $\phi$ of the internal language, we write $(\forall \_ : X) \phi$ for the result
of universally quantifying over a variable that does not appear in $\phi$. Similar 
notation is adopted also for existential quantification. \medskip

We end this section with some examples of Heyting pretoposes equipped 
with families of small maps. Example~\ref{exa:topos} shows that our 
development of internal sheaves applies to elementary 
toposes~\cite[Chapter V]{MacLaneS:shegl}, while Example~\ref{exa:classes}
and Example~\ref{exa:setoid} show that it includes important 
examples for which the ambient category $\catE$ is not an elementary topos.
Note that neither Example~\ref{exa:classes} nor Example~\ref{exa:setoid} 
satisfies the Representability axiom of~\cite[Definition 1.1]{JoyalA:algst}, 
but only the Weak Representabily Axiom, as stated in (S2) above. 
Furthermore, neither of these examples satisfies the additional axiom 
that every object $X$ of $\catE$ has a small diagonal map $\Delta_X : 
X \rightarrow X \times X$.

\begin{example} \label{exa:topos} Consider an elementary topos $\catE$
and let $\mapS$ consist of all maps in $\catE$. It is evident that the
axioms (A1)-(A7) and (P1) are verified, while (S2) is not.
\end{example}

\begin{example} \label{exa:classes}
Consider Constructive Zermelo-Fraenkel set theory ($\mathrm{CZF}$),
presented in~\cite{AczelP:notcst}. We take $\mathcal{E}$ to be 
the exact completion~\cite{CarboniA:somfcr,CarboniA:regec}
of the corresponding category of 
classes~\cite{AwodeyS:preast,GambinoN:premcst}, 
considered as a regular category. By the general theory of exact 
completions, the category $\catE$ is a Heyting pretopos and the 
category of classes of $\mathrm{CZF}$ embeds faithfully in 
it~\cite{CarboniA:somfcr,CarboniA:regec}. 
We write the objects of $\catE$ as $X/r_X$, where $X$ is a 
class and $r_X \subseteq X \times X$ is an equivalence relation on it. 
A map $f : X/r_X \rightarrow A/r_A$ in $\catE$ is a relation 
$f \subseteq X \times A$ that is functional and preserves the 
equivalence relation, in the sense made precise 
in~\cite{CarboniA:somfcr}. We declare a map $f : X/r_X \rightarrow A/r_A$ 
to be small if it fits into a quasi-pullback of the form
\[
\xymatrix{
Y \ar[r]^(.43)k \ar[d]_g & X/r_X \ar[d]^f \\
B \ar@{->>}[r]_(.43)h        & A/r_A }
\]
where $h : B \rightarrow A$ is an epimorphism and $g : Y \rightarrow B$ 
is a function of classes whose fibers are sets. This family satisfies 
the axioms (A1)-(A7) and (S1)-(S2) by a combination of the results on the
category of classes of $\mathrm{CZF}$ in~\cite{AwodeyS:preast,GambinoN:premcst}
with those on small maps in exact completions in~\cite[Section~4]{vandenBergB:exac}. 
\end{example}

\begin{example} \label{exa:setoid}
Consider Martin-L\"of's constructive type theory
with rules for all the standard forms of dependent types and for a 
type universe reflecting 
them~\cite{NordstromB:martt,GambinoN:gentti}.
We take $\catE$ to be the corresponding category of setoids,
which has been shown to be a Heyting pretopos 
in~\cite[Theorem~12.1]{MoerdijkI:typttc}.
We declare a map $f : X \rightarrow A$ in $\catE$ 
to be small if it fits into a quasi-pullback of the form
\[
\xymatrix{
Y \ar[r]^k \ar[d]_g & X \ar[d]^f \\
B \ar@{->>}[r]_h        & A }
\]
where $h : B \rightarrow A$ is an epimorphism and 
$g : Y \rightarrow B$ is a map such that for every $b \in B$
the setoid $g^{-1}(b)$, as defined in~\cite[Section 12]{MoerdijkI:typttc},
is isomorphic to setoid whose carrier and equivalence relation are 
given by elements of the type universe. This family satisfies the 
axioms (A1)-(A7) and (S1)-(S2) by combining the results on display maps 
in~\cite[Section 4]{vandenBergB:exac} 
with those on setoids in~\cite[Section 12]{MoerdijkI:typttc}. 
\end{example}

\section{Lawvere-Tierney sheaves}
\label{sec:myltsheaves}

 Let $\catE$ be a Heyting pretopos
equipped with a family  of small maps $\mapS$ satisfying axioms (A1)-(A7) and
(P1). We define the object of \myemph{small 
truth values} $\myomega$ by letting $\myomega \defeq \smallpower(1)$. 
By the universal property of $\myomega$, it is immediate to see that we 
have a global element $\top : 1 \rightarrow \myomega$. To simplify notation, for $p : \myomega$,
we write $p$ instead of $p = \top$. For example, this 
allows us to write $\lscott p : \myomega \ | \ p \rscott$ instead of 
$\lscott p : \myomega \ | \ p =  \top 
\rscott$. Note, then, that the monomorphism $\lscott p : \myomega \ | \ p \rscott 
\rightarrowtail \myomega$ is the map $\top : 1 \rightarrow \myomega$. 
In a similar fashion, $p \Implies \phi$ is equivalent to  $(\forall \_ \in p) \phi$, for every $p : \myomega$ and every 
formula~$\phi$ of the internal language.  This is because  $(\forall \_ \in p) \phi$ is an abbreviation for
$(\forall x : 1)( x \in p \Rightarrow \phi )$. The internal 
language is used in Definition~\ref{thm:ltcov} to specify what will be
our starting point to introduce a notion of sheaf.

\begin{definition} \label{thm:ltcov} 
Let $(\catE, \mapS)$ be a Heyting pretopos with a family of small maps. 
A \myemph{Lawvere-Tierney coverage} in $\catE$ is a subobject 
$J \rightarrowtail \myomega$ making the following sentences valid in $\catE$
\begin{enumerate}[({C}1)]
\item $J(\top)$,
\item $(\forall p : \myomega)(\forall  q : \myomega) \Big[ \big( 
p  \Rightarrow J(q) \big) \Rightarrow \big( J(p) \Rightarrow J(q) \big)
\Big]$.
\end{enumerate}
\end{definition}

\begin{remark} \label{rem:lawvere-tierney}
Our development of internal sheaves relative to a Lawvere-Tierney
coverage generalises the existing theory of internal sheaves 
relative to a Lawvere-Tierney local operator in an elementary
topos~\cite{MacLaneS:shegl}. Indeed,
when $\catE$ is an elementary topos and $\mapS$ consists of all maps in 
$\catE$, as in Example~\ref{exa:topos}, Lawvere-Tierney coverages are in bijective correspondence with
Lawvere-Tierney local operators~\cite[A.4.4.1]{JohnstoneP:skee}. The
correspondence is given by the universal property of $\myomega$, which is
the subobject classifier of $\catE$, via a pullback square of the form
\[
\xymatrix{
J \ar[r] \ar@{ >->}[d] & 1 \ar[d]^{\top} \\
\myomega \ar[r]_{j} & \myomega }
\]
The verification of the correspondence between the axioms for a 
Lawvere-Tierney coverage and those for a Lawvere-Tierney local
operator is a simple calculation.  Since in 
the general the monomorphism $J \rightarrowtail 
\myomega$ may fail to be small, we focus on Lawvere-Tierney coverages.
\end{remark} 

From now on, we will work with a fixed Lawvere-Tierney coverage as in
Definition~\ref{thm:ltcov}. Our first
step towards defining sheaves is to construct a universal closure operator on
subobjects, that is to say a natural family of functions
\[
\cl_X : \sub(X) \rightarrow \sub(X)\, , 
\] 
for $X \in \catE$, satisfying the familiar monononicity, inflationarity, 
and idempotency properties~\cite[A.4.3]{JohnstoneP:skee}. Note that we do not need to require 
meet-stability, since this follows from the other properties by 
the assumption that the operator is 
natural~\cite[Lemma A.4.3.3]{JohnstoneP:skee}. Naturality of the
operator means that for $f : X \rightarrow Y$, the diagram below commutes
\[
\xymatrix{
\sub(Y) \ar[rr]^{\cl_Y} \ar[d]_{f^*} & & \sub(Y) \ar[d]^{f^*} \\
\sub(X) \ar[rr]_{\cl_X}        & & \sub(X) }
\]
We define the universal closure operator associated to the Lawvere-Tierney
coverage by letting, for $S \rightarrowtail X$
\begin{equation}
\label{equ:cl} 
\cl_X(S) \defeq \big\{ x : X \ | \ 
(\exists p : \myomega) \big( J(p) \land \big( p \Rightarrow S(x) \big) 
\big) \big\} \, . 
\end{equation}

\begin{proposition} \label{thm:uniclos}
The family $\cl_X : \sub(X) \rightarrow \sub(X)$, for 
$X \in \catE$, associated to a Lawvere-Tierney operator is a universal 
closure operator.
\end{proposition} 

\begin{proof} First, we verify that the operator is natural. This is 
immediate, since for a subobject $T \rightarrowtail Y$ we have
\begin{eqnarray*}
\cl_X(f^*T) & = & \big\lscott x : X \ | \ (\exists p : \myomega) 
\big( J(p) \land \big( p \Rightarrow f^*T(x) \big) \big) \big\rscott \\
 & = & f^* \lscott y : Y \ | \  (\exists p : \myomega) \big(
J(p) \land \big( p \Rightarrow T(y) \big) \big)  \big\rscott \\
 & = & f^* ( \cl_Y(T) ) \, . 
\end{eqnarray*}
For inflationarity, let  $S \rightarrowtail X$, $x : X$ and assume
that $S(x)$ holds. Then, define $p : \myomega$ by letting $p \defeq \top$. 
We have that $J(p)$ holds by (C1), and that $p \Rightarrow S(x)$ holds
by assumption.  Monotonicity is immediate by the definition 
in~(\ref{equ:cl}). Idempotence is the only part that is not
straightforward, since it makes use of the Collection Axiom
for small maps. For $S \rightarrowtail X$, we need to show that
$\cl^2(S) \subseteq \cl(S)$. Let $x : X$ and assume that there exists
$p : \myomega$ such that $J(p)$ and $p \Rightarrow \cl_X(S)(x)$ hold.
For $\_ : 1$ and $q : \myomega$, let us define
\[
\phi(\_ \, ,q) \defeq J(q) \land \big( q \Rightarrow S(x) \big) \, . 
\]
By the definition in~(\ref{equ:cl}), $p \Rightarrow \cl_X(S)(x)$ implies
\[
(\forall \_ \in p)(\exists q : \myomega) \phi(\_ \, , q) \, . 
\]
We can apply Collection and derive the existence of $u : 
\smallpower(\myomega)$
such that 
\[
(\forall \_ \in p)(\exists q \in u) \phi(\_ \, ,q)
\land 
(\forall q \in u)(\exists \_ \in p) \phi(\_ \, ,q) \, . 
\]
Define $r : \myomega$ by $r \defeq \bigcup u$. We wish to show that
$J(r)$ and $r \Rightarrow S(x)$ hold, which will allow us
to conclude $\cl S_X(x)$, as required. To prove that $J(r)$ holds,
we observe that
\begin{eqnarray*}
p & \Rightarrow  & (\exists q \in u) \, J(q)  \\
  & \Rightarrow  & (\exists q : \myomega) \big( J(q) \land q \subseteq r \big) \\
  & \Rightarrow & J(r) \, . 
\end{eqnarray*}
Therefore $p \Rightarrow J(r)$. Since we have $J(p)$ by hypothesis,
Axiom~(C2) for a Lawvere-Tierney coverage allows to derive $J(r)$,
as required. By definition of $r : \myomega$, we note that $r \Rightarrow 
S(x)$ holds if and only if for every $q \in u$ we have $q \Rightarrow S(x)$.
But since $q \in u$ implies $q \Rightarrow S(x)$, we obtain
$r \Rightarrow S(x)$, as desired.
\end{proof} 

Let us point out that the proof of Proposition~\ref{thm:uniclos} makes use of the
Collection Axiom for small maps. This is analogous to the application of the 
type-theoretic collection axiom in~\cite[Lemma 5.4]{GambinoN:gentti}.

\begin{remark} \label{thm:closvscov}
Given a universal closure operator on $\catE$, it is possible to 
define a Lawvere-Tierney coverage $J \rightarrowtail \myomega$ by taking
$J$ to be the closure of the subobject $\{ p : \myomega \ | \ p \}
\rightarrowtail \myomega$. The closure operator
induced by $J$ coincides with the given one
if and only if the latter satisfies the equation in~(\ref{equ:cl}). 
Therefore, we are considering here only a special class of universal closure
operators. Restricting to this special class allows us to develop 
a treatment of sheaves  without 
assuming additional axioms for small maps, such as that asserting that 
every monomorphism is small. 
Focusing on the class of universal closure operations
determined by Lawvere-Tierney
coverages captures an appropriate level of generality. First, 
as we will see below, they 
correspond precisely to the Grothendieck sites considered 
in~\cite{GambinoN:asssft}. Secondly, as
we will see in Section~\ref{sec:asssft} and Section~\ref{sec:smams},
they allow us to develop the construction of internal sheaves.

For an example of a universal closure operator that does not seem to 
satisfy the equation in~(\ref{equ:cl}), consider the double-negation 
universal closure 
operator, defined by mapping  a subobject $S \rightarrowtail X$ into the 
subobject $\{ x : X \ | \ \neg \neg S(x) \} \rightarrowtail X$. The induced 
Lawvere-Tierney coverage is given by $\{ p : \Omega \ | \ \neg \neg p \} 
\rightarrowtail \Omega$. 
Without the assumption of further axioms for small maps, 
the universal closure operator associated to this Lawvere-Tierney coverage 
may be different from the double-negation closure operator. Within our
setting, therefore, the category-theoretic counterpart of the 
double-negation translation has to be developed  by 
working with the universal closure operator
associated to the Lawvere-Tierney coverage $\{ p : \Omega \ | \ \neg \neg p \} \rightarrowtail \Omega$.
\end{remark} 

We shall be particularly interested in the closure of the membership
relation $\ni_X \rightarrowtail \smallpower(X) \times X \rightarrow X$,
which we are going to write as 
\[
\nibar[X] \rightarrowtail \smallpower(X) \times X \rightarrow 
\smallpower(X) \, . 
\]
For $x : X$ and $s : \smallpower(X)$, the definition in~(\ref{equ:cl}) 
implies
\[
x \inbar s \Iff (\exists p : \myomega) \big( J(p) \land 
\big( p \Rightarrow 
x \in s \big) \big) \, ,
\]

Having defined a universal closure operator on $\catE$, we can define
a notion of sheaf and introduce the family of maps that we will consider
as small maps between sheaves. In order to do this, we define a monomorphism 
$m : B \rightarrowtail A$ to be \emph{dense} if it holds that 
$\cl_A(B) = A$. 

\begin{definition}  \label{thm:ltsheaf}
An object $X$ of $\catE$ is said to be a 
\myemph{separated} if  for every dense monomorphism 
$m : B \rightarrowtail A$ the function
\begin{equation}
\label{equ:shemap}
\hom(m,X) : \hom(A,X) \rightarrow \hom(B,X) \, , 
\end{equation}
induced by composition with $m$, is injective. Equivalently,
$X$ is a separated if and only if for every map $v : B \rightarrow X$
there exists at most one extension $u : A \rightarrow X$ making the 
following diagram commute
\begin{equation}
\label{equ:diagshe}
\xymatrix{
B \ar[r]^(.45){v} \ar@{ >->}[d]_{m} & X \\
A \ar@/_/[ur]_{u} & }
\end{equation}
We say that $X$ is a \myemph{sheaf} if the 
function in~(\ref{equ:shemap}) is bijective. Equivalently,~$X$ is a sheaf if and only if
every map $v : B \rightarrow X$ has a unique 
extension $u : A \rightarrow X$ as in~(\ref{equ:diagshe}).
\end{definition}

Definition~\ref{thm:jsmall} defines what it means for a map in $\catE$ 
to be \myemph{locally small}. As explained further in 
Section~\ref{sec:class}, the idea 
underlying this notion is that a map is locally small if and only
if each of its fibers contains a small dense subobject. 
By their very definition, locally small maps are stable under 
pullback, since the properties of the defining diagrams all are. 
 
\begin{definition} \label{thm:jsmall}
A map $f : X \rightarrow A$ in $\catE$ is said to be \emph{locally small} 
if its graph $\Gamma(f) : X \rightarrowtail A \times X \rightarrow A$ 
fits in a diagram of the form
\[
\xymatrix{
T \ar[r] \ar@{ >->}[d] & X \ar@{ >->}[d] \\
B \times X \ar[r] \ar[d]_{\pi_B} & A \times X \ar[d]^{\pi_A} \\
B \ar@{->>}[r]_h & A }
\]
where $h : B \twoheadrightarrow A$ is an epimorphism,
$T \rightarrowtail B \times X \rightarrow B$ is
an indexed family of small subobjects of $X$, and
the canonical monomorphism $T \rightarrowtail B \times_A X$
is dense. An indexed family of subobjects $S \rightarrowtail 
A \times X \rightarrow X$ is said to be an 
\myemph{indexed family of locally small subobjects} if the 
composite map $S \rightarrow A$ is locally small.  
\end{definition}

We write $\JSh$ for the full subcategory of $\catE$ whose objects are
sheaves, and~$\mapS_J$ for the family of locally small maps in $\JSh$. 
The aim of the remainder of the paper is to prove the following result.
As fixed in Section~\ref{sec:asthp}, a family of small maps is
required to satisfy only axioms~(A1)-(A7) and (P1). 

\begin{theorem} \label{thm:main} 
Let $(\catE,\mapS)$ be a Heyting pretopos with a 
family of small maps. For every Lawvere-Tierney coverage 
$J$ in $\catE$, $(\JSh, \mapS_J)$ is a Heyting pretopos equipped with
a family of small maps. Furthermore, if $\mapS$ satisfies the
Exponentiability and Weak Representability axioms, so does $\mapS_J$. 
\end{theorem}

Theorem~\ref{thm:main} concerns only the axioms for small maps discussed in Section~\ref{sec:asthp}. To obtain
sheaf models for constructive, intuitionistic, and classical set theories, it is necessary to consider additional axioms
for small maps in~$\catE$ and show that they are preserved by the construction of internal sheaves. We expect that
the preservation of some axioms for small maps by the construction of internal sheaves requires the assumption of additional conditions on the Lawvere-Tierney coverage. We leave the treatment of these issues to future research.  
Let us point out  that our preference for assuming the Weak Representability axiom of~\cite{vandenBergB:uniaas} rather than the Representability axiom of~\cite{JoyalA:algst} is due to the fact that it does not seem possible to prove that the latter is preserved by our construction without the assumption of additional axioms for
small maps.  \medskip

Before developing the theory required to prove
Theorem~\ref{thm:main}, we explain how this result subsumes the
treatment of sheaves for a Grothendieck site. 
Let $\mathbb{C}$ be a small internal category in $\catE$.
Smallness of $\mathbb{C}$ means that both its objects $\mathbb{C}_0$ 
and its arrows $\mathbb{C}_1$ are given by 
small objects in $\catE$. We write category 
$\Psh[\catE](\mathbb{C})$ of internal presheaves over~$\mathbb{C}$.
It is well-known that $\Psh[\catE](\mathbb{C})$ is a Heyting pretopos. 
When $\mathbb{C}_0$ and $\mathbb{C}_1$ have small diagonals, so that
both the equality of objects and that of arrows is given by a small
monomorphism, it is possible to equip $\Psh[\catE](\mathbb{C})$
with a family of small maps, consisting of the internal natural
transformations that are pointwise small maps in 
$\catE$~\cite{MoerdijkI:typttc,WarrenM:catc}. 
Lawvere-Tierney coverages in $\Psh[\catE](\mathbb{C})$ are in
bijective correspondence with Grothendieck coverages with small
covers on $\mathbb{C}$, as defined in~\cite{GambinoN:asssft}. To 
explain this correspondence, we need to recall some terminology and 
notation. A \myemph{sieve}~$P$ on $a \in \mathbb{C}$ is 
a subobject $P \rightarrow \yon(a)$, where we write $\yon(a)$ for the Yoneda
embedding of $a$. Such a sieve can be identified 
with a family of 
arrows with codomain $a$ that is closed under composition, in the sense 
that for every pair of composable maps $\phi : b \rightarrow a$ and 
$\psi : c \rightarrow b$ in $\mathbb{C}$, $\phi \in P$ implies 
$\phi \, \psi \in P$. For a  sieve $P$ on $a$ and 
arrow $\phi : b \rightarrow a$, we write $P \cdot \phi$ for 
the sieve on $b$ defined by letting
\begin{equation}
\label{equ:res} 
P \cdot \phi \defeq 
\{ \psi : c \rightarrow b \ | \ \phi \, \psi \in P \} \, .
\end{equation}
Recall from~\cite{GambinoN:asssft} that a \myemph{Gro\-then\-dieck coverage with small covers} 
on $\mathbb{C}$ consists of a family $(\cov(a) \ | \ a \in \mathbb{C})$ such 
that elements of $\cov(a)$ are small sieves, and the conditions
of Maximality (M), Local Character (L), and Transitivity (T) hold:
\begin{enumerate}
\item[(M)] $M_a \in \cov(a)$.
\item[(L)] If $\phi : b \rightarrow a$ and $S \in \cov(a)$, then $S \cdot \phi
\in \cov(b)$. 
\item[(T)] If $S \in \cov(a)$, $T$ is a small sieve on $a$, and for all
$\phi : b \rightarrow a \in S$ we have $T \cdot \phi \in \cov(b)$, then 
$T \in \cov(a)$.
\end{enumerate}
The object $\myomega$ in $\Psh[\catE](\mathbb{C})$ is given by 
\[
 \myomega(a) \defeq \lscott  S \rightarrowtail \yon(a) \ 
| \ S \text{ small} \rscott \, . 
\]
Therefore, a Grothendieck coverage with small covers can be identified with 
a family of subobjects $\cov(a) \rightarrowtail \myomega(a)$. Condition 
(L) means that this family is a subpresheaf of $\myomega$, while conditions
(M) and (T) for a Grothendieck coverage are the rewriting of 
conditions (C1) and (C2) for a Lawvere-Tierney coverage. By instanciating
the general notion of Definition~\ref{thm:ltsheaf}  we obtain a notion 
of sheaf, which can be shown to be equivalent to the familiar notion of
a sheaf for a Grothendieck coverage, and a notion of small map. Writing
$\Sh[\catE](\mathbb{C}, \cov)$ for the category of internal sheaves,
and $\mapS(\mathbb{C}, \cov)$ for the corresponding family of small
maps, Theorem~\ref{thm:main} implies the following result.

\begin{corollary} \label{thm:cormain} 
Let $(\catE, \mapS)$ be a Heyting pretopos with
a family of small maps. Let $(\mathbb{C}, \cov)$ be a small category with small
diagonals equipped with a Gro\-then\-dieck coverage with small covers. 
Then $\big( \Sh[\catE](\mathbb{C}, \cov), \mapS(\mathbb{C}, \cov) \big)$ 
is a Heyting pretopos with a family of small maps. Furthermore, if $\mapS$ satisfies the
Exponentiability and Weak Representability axioms, so does $\mapS_J$. 
\end{corollary} 

We conclude the section by explaining in more detail how Corollary~\ref{thm:cormain} relates to
similar results in the existing literature. First, it extends the main result of~\cite{GambinoN:asssft}
by considering the preservation not only of the structure of a Heying pretopos, but also of the
axioms for small maps. Corollary~\ref{thm:cormain} neither implies nor is implied by 
the corresponding results in~\cite{vandenBergB:shept,MoerdijkI:typttc}. This is because
the axiomatisation of small maps considered here is rather different to the one
considered in~\cite{vandenBergB:shept,MoerdijkI:typttc}. The axiomatisation in~\cite{vandenBergB:shept,MoerdijkI:typttc} takes as motivating example
a class of small maps within the category of setoids  that is different from the one considered
here. Furthermore, locally small maps  and  pointwise small maps between sheaves coincide
in~\cite{vandenBergB:shept,MoerdijkI:typttc}, while they do not here.
Finally, the result on construction of internal sheaves announced in~\cite{vandenBergB:uniaas},
which considers the preservation of axioms for small maps that we do not study here, involves 
additional assumption on the Grothendieck sites, which are not part of the hypotheses of 
Corollary~~\ref{thm:cormain}.

\section{Classification of locally small subobjects}
\label{sec:class}

We begin by characterising locally small maps in the
internal language, analogously to how small maps are 
characterised in~(\ref{equ:logchasmall}). For each object $X$, define 
an equivalence relation $R \rightarrowtail \smallpower(X) \times 
\smallpower(X)$ by letting
\[
R \defeq \lscott (s,t) : \smallpower(X) \times \smallpower(X) \ | \ 
(\forall x:X) \big( x \inbar s \iff x \inbar t  \big) \rscott \, .
\]
Informally, $R(s,t)$ holds whenever $s$ and $t$ have the same closure.
Using the exactness of the Heyting pretopos $\catE$, we define $\peejay(X)$ 
as the quotient of $\smallpower(X)$ by $R$, fitting into an exact diagram of
the form
\begin{equation}
\label{equ:relequ}
\xymatrix{
R \ar@<0.75ex>[r]^(.42){\pi_1} \ar@<-0.75ex>[r]_(.42){\pi_2} & 
\smallpower(X) \ar@{->>}[r]^(.44){[\; \cdot \; ]} & 
\peejay(X) \, . }
\end{equation}
The quotient map $\smallpower(X) \twoheadrightarrow 
\peejay(X)$ is to be 
interpreted as performing the closure of a small subobject of $X$. 

\begin{remark} The exactness of the Heyting pretopos $\catE$
is exploited here in a crucial way to define $\peejay(X)$. In
particular, without further assumptions on the Lawvere-Tierney
coverage, the equivalence relation in~(\ref{equ:relequ}) cannot
be shown to be given by a small monomorphism.
\end{remark}

We define a new indexed family of subobjects of $X$,
$\nij \rightarrowtail \peejay(X) \times X \rightarrow \peejay(X)$, by 
letting, for $x : X$ and $p : \peejay(X)$, 
\[
x \inj p \iff (\exists s : \smallpower X)  
\big( p = [s] \land x \inbar s \big) \, . 
\]
In particular, for $x:X$ and $s :\smallpower X$, this gives $x \inj [s] 
\iff  x \inbar s$, which is to say that the following squares are pullbacks
\begin{equation}
\label{equ:pbpj}
\xymatrix{
\nibar \ar[r] \ar@{ >->}[d] & \nij \ar@{ >->}[d] \\
\smallpower(X) \times X \ar[r] \ar[d] & \peejay(X) \times X \ar[d] \\
\smallpower(X) \ar[r]_{[\, \cdot \, ]} & \peejay(X) }
\end{equation} 
This definition of $\inj$ and of the relation $R$ imply that
$\peejay(X)$ satisfies a form of extensionality, in the sense
that for $p, q : \peejay(X)$ it holds that
\[
p = q  \Iff (\forall x:X) ( x \inj p \Iff x \inj q) \, . 
\]
From diagram (\ref{equ:pbpj}), we also see that $\inj$ is closed in 
$X \times \peejay X$, since it is closed when pulled back along an 
epimorphism. Given a map $f : X \rightarrow A$, it is convenient to define,
for $a : A$, $s : \smallpower(X)$,
\[ 
s \approx f^{-1}(a) \defeq
(\forall x : X) 
\left[ 
(x \in s \Implies f(x) = a) \land  
(f(x) = a \Implies x \inbar s) \right] \, . 
\]

\begin{lemma} \label{thm:logjsmall}
A map $f : X \rightarrow A$ is locally small if and only if the following 
sentence is valid: 
\begin{equation}
\label{equ:logsmall}
(\forall a : A) 
(\exists s : \smallpower X) s \approx f^{-1}(a) \, . 
\end{equation}
\end{lemma} 

\begin{proof} First, let us assume that $f : X \rightarrow A$ is 
locally small. By Definition~\ref{thm:jsmall} there is a diagram
\[
\xymatrix{
T \ar[r] \ar@{ >->}[d] & X \ar@{ >->}[d] \\
B \times X \ar[r] \ar[d] & A \times X \ar[d] \\
B \ar@{->>}[r]_h & A } 
\]
where $T \rightarrowtail B \times X \rightarrow B$ is an  indexed family
of small subobjects of $X$,  $h : B \twoheadrightarrow A$ is an epimorphism, and
the monomorphism $T \rightarrowtail B \times_A X$ is dense.
There is a classifying map $\chi_T : B \rightarrow \smallpower(X)$ such that
$x \in \chi_T(b) \iff T(b,x)$, and the commutativity of the diagram implies 
\[
x \in \chi_T(b) \Implies f(x) = h(b)  \, ,
\]
while the density of  $T \rightarrowtail B \times_A X$ implies 
\[
f(x) = h(b) \Rightarrow x \inbar \chi_T(b) \, . 
\]
Given $a : A$, there exists $b : B$ such that $h(b) = a$ and so, 
defining $s \defeq \chi_T(b)$, we obtain the data required
to prove the statement. For the converse 
implication, assume~(\ref{equ:logsmall}).
We define
\[
B \defeq \lscott (a,s) : A \times \smallpower(X) \ | \ 
 s \approx f^{-1}(a) \rscott 
\]
and 
\[
T \defeq \lscott \big( (a,s), x \big) : B \times X \ | \ x \in s \} \, . 
\]
We have that $T \rightarrowtail B \times X \rightarrow B$ is small by 
construction and that the projection $h : B \rightarrow A$ is an 
epimorphism by hypothesis. Since
\[
B \times_A X \iso \lscott \big( (a,s), x \big) : B \times X \ |  \ f(x) = a 
\rscott \, , 
\]
we have that $T \rightarrowtail B \times_A X$ is dense by the definition 
of $B$. 
\end{proof}

\begin{remark}  We could have
considered maps $f:X \rightarrow A$ satisfying the condition
\[
(\forall a : A)(\exists s: \smallpower X)(\forall x:X.) 
f(x) = a \Iff x \inbar s \, . 
\]
This amounts to saying that each fiber of $f$ is the closure of a 
small subobject. Let us call such maps \myemph{closed-small}. For maps with codomain a separtated object, and so for maps
with codomain a sheaf, these
definitions coincide, so either would give our desired class of small maps 
on $\JSh$. However, considered on the whole of $\catE$, they 
generally give different classes of maps, each retaining different 
properties. For instance, closed-small subobjects of $X$ are classified 
by $\peejay(X)$, while locally small subobjects may not be classified. On 
the other hand, locally small maps satisfy axioms (A1)-(A7) in 
$\catE$, whereas the closed-small 
maps may not, since the identity map on a non-separated object 
will not be closed-small. Thus, the choice of either as the extension of 
local smallness from $\JSh$ to $\catE$ is a matter of convenience.
 \end{remark}

We conclude this section by showing how locally
small closed subobjects can be classified. 
This is needed in Section~\ref{sec:peejay-and-asf} 
for the proof of the associated sheaf functor theorem.

\begin{lemma} \label{thm:usefullemma}
An indexed family of subobjects 
$S \rightarrowtail A \times X \rightarrow A$ 
is an indexed family of locally small subobjects if and only if there exists
a diagram of the form
\[
\xymatrix{
T \ar[r] \ar@{ >->}[d] & S \ar@{ >->}[d] \\
B \times X \ar[r] \ar[d]_{\pi_B} & A \times X \ar[d]^{\pi_A} \\
B \ar@{->>}[r]_h & A }
\]
where $h : B \twoheadrightarrow A$ is an epimorphism,
$T \rightarrowtail B \times X \rightarrow S$ is
an indexed family of small subobjects of $X$, and
the canonical monomorphism $T \rightarrowtail B \times_A S$
is dense.  
\end{lemma}

\begin{proof} Assume $S \rightarrowtail A \times X \rightarrow A$ 
to be an indexed family of locally small subobjects. 
By Definition~\ref{thm:jsmall}, this means that the composite
$S \rightarrow A$ is a locally small map, which in turn implies
that there exists a diagram of the form 
\begin{equation}
\label{equ:diaggg}
\xymatrix{
T \ar[r] \ar@{ >->}[d] & S \ar@{ >->}[d] \\
B \times S \ar[r] \ar[d]_{\pi_B} & A \times S \ar[d]^{\pi_A} \\
B \ar@{->>}[r]_h & A }
\end{equation}
where $h : B \twoheadrightarrow A$ is an epimorphism,
$T \rightarrowtail B \times S \rightarrow S$ is
an indexed family of small subobjects of $S$, and
the canonical monomorphism $T \rightarrowtail B \times_A S$
is dense. The composite 
\[
\xymatrix{
T \ar@{ >->}[r] 
& B \times S \ar@{ >->}[r] & B \times A \times X \ar[r] 
& B \times X }
\]
can be shown to be a monomorphism using the commutativity of
the diagram in~(\ref{equ:diaggg}). We obtain an 
indexed family of small subobjects $T \rightarrowtail B 
\times X \rightarrow B$, which clearly satisfies the required
property. 
\end{proof}

Observe 
that $\nij \rightarrow \peejay(X) \times X \rightarrow \peejay(X)$
is a family of locally small subobjects of $X$, as witnessed by the diagram
\[
\xymatrix{
\ni \ar[r] \ar@{ >->}[d] & \nij \ar@{ >->}[d] \\
\smallpower(X) \times X \ar[r] \ar[d] & \peejay(X) \times X \ar[d] \\
\smallpower(X) \ar@{->>}[r]_{[\, \cdot \, ]} & \peejay(X) }
\]
The pullback of $\nij \rightarrowtail \peejay(X) \times 
X \rightarrow \peejay(X)$ along a map~$\chi : A \rightarrow \peejay(X)$ 
is therefore an 
$A$-indexed family of locally small closed subobjects of $X$. We 
write $\ssubj(X)(A)$ for the lattice of such subobjects.

\begin{proposition} \label{what-peejay-classifies} For every object
$X$, the object $\peejay(X)$ classifies indexed families of locally small 
closed subobjects of $X$, which is to say that for every such
family $S \rightarrowtail A \times X \rightarrow A$ there exists
a unique map $\phi_S : A \rightarrow \peejay(X)$ such that both
squares in the diagram
\[
\xymatrix{
S \ar[rr] \ar@{ >->}[d] & & \nij \ar@{ >->}[d] \\
A \times X \ar[rr]^(.45){\phi_S \times 1_X} \ar[d]_{\pi_A} & & \peejay(X) 
\times X \ar[d]^{\pi_{\peejay(X)}} \\
A \ar[rr]_{\phi_S} & & \peejay(X) }
\]
are pullbacks. Equivalently, the functions
\[
\hom(A, \peejay(X)) \rightarrow \ssubj(X)(A)
\]
given by pulling back $\nij \rightarrowtail \peejay(X) \times X \rightarrow
\peejay(X)$ are a family of bijections, natural in $A$. 
\end{proposition}

\begin{proof} Let an indexed family as in the statement be given,
together with data as in Lemma~\ref{thm:usefullemma}. By the universal 
property of $\smallpower(X)$, we get a classifying map $\chi_T : B \rightarrow 
\smallpower(X)$. Using the naturality of the closure operation, we
obtain
\[
(\chi_T \times 1_X)^* \cl(\ni_X)  = 
\cl\big( (\chi_T \times 1_X)^*(\ni_X) \big) =
\cl(T) =
B \times_A S \, ,
\]
where the last equality is a consequence of the density of $T$ in
$B \times_A S$. Therefore, we have a sequence of pullbacks of the
form
\[
\xymatrix{
\nij \ar@{ >->}[d] & \nibar \ar[l] \ar@{ >->}[d] & B \times_A S \ar[r] \ar[l] 
\ar@{ >->}[d]  & S \ar@{ >->}[d] \\
\peejay(X) \times X \ar[d] & \smallpower(X) \times X  \ar[l] \ar[d] &
 B \times X \ar[r]  \ar[l] \ar[d] & A \times X \ar[d] \\
\peejay(X) & \smallpower(X) \ar@{->>}[l]^{[\, \cdot \, ]}  & B 
\ar@{->>}[r]_h \ar[l]^{\chi_T} & A }
\]
This, combined  with the epimorphism $h : B 
\rightarrow A$, allows us to show that
\[
(\forall a : A)
(\exists p : \peejay(X))  
(\forall x : X) \, \big(
S(a,x) \Iff x \inj p \big) \, . 
\]
But the definition of $\peejay(X)$ as a quotient ensures that for a given $a : A$,
there is a unique $p : \peejay(X)$ satisfying $S(a,x) \Iff x \inj p$
for all $x : X$. Hence, by functional completeness
we obtain the existence of a map 
$\phi_S : A \rightarrow \peejay(X)$ such that
\[
(\forall a : A)
(\forall x : X) \, 
S(a,x) \Iff x \inj \phi_S(x) \, , 
\]
as required.
\end{proof} 

\section{The associated sheaf functor theorem} 
\label{sec:peejay-and-asf}
\label{sec:asssft}

We are now ready to define the 
associated sheaf functor $\asf : \catE \rightarrow \JSh$,
the left adjoint to the inclusion 
$\incl : \JSh \rightarrow \catE$ of the full
subcategory of sheaves into $\catE$. 
Given  $X \in \catE$, define $\sigma : X \rightarrow \peejay(X)$ to be the 
composite
\[
\xymatrix{
X \ar[r]^(.39){ \{ \, \cdot \, \} }  & \smallpower(X) 
\ar@{->>}[r]^(.45){ [ \, \cdot \, ] }
& \peejay(X) \, . }
\]
To define the associated sheaf functor, first factor $\sigma : X 
\rightarrow \peejay(X)$ as an epimorphism $X \twoheadrightarrow X'$ 
followed by a monomorphism $X' \rightarrowtail \peejay(X)$. Then,
define $\asf(X)$ to be the closure of the subobject $X'$ in $\peejay(X)$
\[
         \asf(X) \defeq \cl(X') \, . 
\]
The unit of the adjunction $\eta_X : X \rightarrow \asf(X)$ is then 
defined as the composite of $X \twoheadrightarrow X'$ with the inclusion 
$X' \rightarrowtail \cl(X')$. We need to show that~$\asf(X)$ is indeed a 
sheaf and that it satisfies the appropriate universal property. The 
proof of the former involves 
the verification that $\peejay(X)$ is a sheaf. This, in turn, requires 
further analysis of the notion of locally small map, which we carry out 
next.

We say that a map is \myemph{dense} if it factors as an epimorphism 
followed 
by a dense monomorphism. For monomorphisms this definition 
agrees with the definition of dense mo\-no\-mor\-phism given in 
Section~\ref{sec:myltsheaves}. It is 
immediate to see that the pullback of a dense map is again dense, 
and a direct calculation shows that dense maps are closed under 
composition. It will be convenient to introduce some additional 
terminology: we refer to a commutative square of the form 
\[
\xymatrix{
Y \ar[d] \ar[r] & X \ar[d] \\
B \ar[r] & A }
\]
such that the canonical map $Y \rightarrow B \times_A X$ is dense
as a \emph{local quasi-pullback}. Note that every dense map $h : 
B \rightarrow A$ fits into a local quasi-pullback of the form
\[
\xymatrix{
B \ar[r]^h \ar@{ >->}[d]_m & A \ar[d]^{1_A} \\
B' \ar@{->>}[r]_p & A }
\]
where $m : B \rightarrowtail B'$ and $p : B' \twoheadrightarrow A$ 
are respectively the dense monomorphism and the 
epimorphism forming the factorisation of $h$. The diagram is a 
local quasi-pullback 
because the map $B \rightarrow B' \times_A A$ is $m$ itself. Let us also
observe that if  both squares in a diagram
\[
\xymatrix{
Z \ar[r] \ar[d] & Y \ar[r] \ar[d] & X \ar[d] \\
C \ar[r]        & B \ar[r]        & A }
\]
are local quasi-pullbacks, then the whole rectangle is also a local
quasi-pullback. This implies
that any finite pasting of local quasi-pullbacks is again a 
local quasi-pullback. We establish a very useful factorisation for dense 
mononomorphisms.

\begin{lemma} \label{lemma:handy-factorisation} Every dense monomorphism 
$h : B \rightarrowtail A$ can be factored as 
\[
\xymatrix{
B \ar@{ >->}[rr]^m \ar@{ >->}[dr]_h & & B' \ar@{->>}[dl]^p  \\
            & A & }
\] 
where $m : B \rightarrowtail B'$ is a small dense monomorphism 
and $p : B' \twoheadrightarrow A$ is an epimorphism. 
\end{lemma}

\begin{proof} Let us define $B' \defeq \lscott 
(p,a) : \myomega \times A \ | \ p \Rightarrow B(a) \rscott$.
By the definition of closure in~(\ref{equ:cl}) and
the assumption that $h$ is dense, the projection $\pi_2 : B'
\rightarrow A$ is an epimorphism. Furthermore, there
is a monomorphism $m : B \rightarrow B'$ defined by mapping
$b : B $ into $(\top, b) : B'$. Diagrammatically, we
have a pullback of the form
\[
\xymatrix{
B \ar[r] \ar@{ >->}[d]_m   & 1 \ar@{ >->}[d]^{\top} \\
B' \ar[r]_{\pi_1}  &  J  }
\]
The map $\top : 1 \rightarrow J$ is small and
dense. It is small because it is the pullback
of the map $\top : 1 \rightarrow \myomega$, which
is small by the definition of $\myomega$, 
along the inclusion $J \rightarrowtail \myomega$.
It is dense by the very definition of closure in~(\ref{equ:cl}).
Therefore, preservation of smallness and density along pullbacks 
implies that $m : B \rightarrow B'$ is small and dense, as required.   
\end{proof} 

\begin{remark} Lemma~\ref{lemma:handy-factorisation} exploits in
a crucial way the fact that the closure operation, and hence the
notion of density, are determined by a Lawvere-Tierney coverage.
Indeed, it does not seem possible to prove an analogue of its 
statement for arbitary closure operations without assuming 
further axioms for small maps. 
\end{remark}

\begin{lemma} \label{j-small-as-j-q-pb} A map $f:X \rightarrow A$ is 
locally  small if and only if it fits into a local quasi-pullback 
square of the form
\begin{equation}
\label{equ:givendiag}
\xymatrix{
Y \ar[r]^k \ar[d]_g & X \ar[d]^{f} \\
B \ar@{->>}[r]_h        & A }
\end{equation}
where $h : B \twoheadrightarrow A$ is an epimorphism and $g : Y \rightarrow B$
is small. 
\end{lemma}

\begin{proof} 
Assuming that $f:X \rightarrow A$ is locally small,
the required diagram is already given by Definition~\ref{thm:jsmall}. 
For the converse implication, consider
a diagram as in~(\ref{equ:givendiag}). Consider the factorisation of
the canonical map $Y \rightarrow B \times_A X$ as an epimorphism followed
by a monomorphism, say $Y \twoheadrightarrow T \rightarrowtail B \times_A X$.
We have a diagram of form
\[
\xymatrix{
Y \ar@{->>}[r] \ar@{ >->}[d] & T \ar@{ >->}[d] \ar[r] & X \ar@{ >->}[d] \\
B \times Y \ar[r] \ar[d]_{\pi_B} & B \times X \ar[d]^{\pi_B} \ar[d]
\ar[r]  & 
A \times X \ar[d]^{\pi_A} \\
B \ar[r]_{1_B} & B \ar@{->>}[r]_h & A }
\]
By the Quotients Axiom for small maps 
$T \rightarrowtail B \times X \rightarrow B$ is a family of
small subobjects. Furthermore, 
$T \rightarrowtail B \times_A X$ is dense by construction. 
\end{proof}

Note that the map $k : Y \rightarrow X$ in Lemma~\ref{j-small-as-j-q-pb} is 
also dense, since it is the composition of the pullback of an epimorphism 
with a dense map. 

\begin{proposition} \label{lemma:J-smallness-is-local} 
 \label{lemma:J-smallness-is-local-2}
 \label{j-small-as-j-q-pb-2}
A map $f : X \rightarrow A$
is locally small if and only if it fits into a pullback diagram of the form
\begin{equation}
\label{equ:testproof}
\xymatrix{
Y \ar[r]^k \ar[d]_g & X \ar[d]^f \\
B \ar[r]_h  & A }
\end{equation}
where $h : B \rightarrow A$ is dense and $g : Y \rightarrow B$ is locally
small.
\end{proposition}

\begin{proof} Assume to be given a diagram as (\ref{equ:testproof}).
By the factorisation of dense maps as
epimorphisms followed by dense monomorphisms  and 
Lemma~\ref{j-small-as-j-q-pb}, it is sufficient to prove the 
statement when $h : B \rightarrow A$ is a dense monomorphism. We construct
a diagram of the form
\[
\xymatrix{
B' \times_C Z    \ar[r] \ar[d] \ar@{}[dr]|{(5)} & 
Z                \ar[r] \ar[d] \ar@{}[dr]|{(3)} & 
Y                \ar[r]^k \ar[d]_g \ar@{}[dr]|{(1)} & 
X \ar[d]^f                         \\
B'               \ar[r] \ar[d]_{n} \ar@{}[drr]|{(4)} & 
C                \ar@{->>}[r]        & 
B                \ar[r]^h \ar[d]_m  \ar@{}[dr]|{(2)} & 
A                \ar[d]^{1_A} \\ 
A''              \ar@{->>}[rr]_{p'}       &     
                                    & 
A'                \ar@{->>}[r]_{p}        & 
A }
\]
The given diagram is in (1). 
First, factor $h : B \rightarrow A$ using Lemma~\ref{lemma:handy-factorisation}
as a small dense monomorphism $m : B \rightarrowtail A'$ followed by an 
epimorphism $p : A' \rightarrow A$. The resulting commutative square in (2) 
is a local quasi-pullback since $m$ is dense. Next, apply 
Lemma~\ref{j-small-as-j-q-pb} to
the locally small map $g : Y \rightarrow B$ so as to obtain the
local quasi-pullback in~(3) with $Z \rightarrow C$ a small map and 
$C \twoheadrightarrow B$
an epimorphism. Next, we apply the Collection Axiom to the small map
$m : B \rightarrow A'$ and the epimorphism $C \rightarrow B$, so as to obtain
the quasi-pullback in~(4) where
$n : B' \rightarrow A''$ is a small map and $p' : A'' 
\twoheadrightarrow A'$ is an epimorphism.
Finally, we construct the pullback square in~(5)
to obtain  another small map $B' \times_C Z \rightarrow B'$. 
The whole diagram is a local quasi-pullback, and we can
apply Lemma~\ref{j-small-as-j-q-pb} to deduce that
$f : X \rightarrow A$ is locally small. The converse implication
is immediate since locally small maps are stable under pullback. 
\end{proof} 

Corollary~\ref{thm:addthis} shows that the assumption of
being an epimorphism for the map $h : B \rightarrow A$ in
Lemma~\ref{j-small-as-j-q-pb} can be weakened. 

\begin{corollary} \label{thm:addthis}
If a map $f : X \rightarrow A$ fits into a local quasi-pullback diagram 
of the form
\begin{equation}
\label{equ:secondgivendiag}
\xymatrix{
Y \ar[d]_g \ar[r]^k & X \ar[d]^f \\
B \ar[r]_h            & A }
\end{equation} 
where $g : Y \rightarrow B$ is small and $h : B \rightarrow A$ is
dense, then $f : X \rightarrow A$ is locally small. 
\end{corollary}

\begin{proof} Given the diagram in~(\ref{equ:secondgivendiag}), 
we construct the following one
\[
\xymatrix{
Y \ar[r] \ar[d]_{g} & B \times_A X \ar[r] \ar[d]^{g'} & X 
\ar[d]^f \\
B \ar[r]_{1_B}        & B            \ar[r]_h        & A }
\]
The left-hand side square is a local quasi-pullback by assumption.
So, by Lemma~\ref{j-small-as-j-q-pb}, the map $g' : B \times_A X \rightarrow
B$ is locally small. The right-hand side square is a pullback by
definition. So, by Proposition~\ref{lemma:J-smallness-is-local}, 
the map $f : X \rightarrow A$ is locally small, as desired.
\end{proof}

We exploit our characterisation of locally small maps in the 
proof of the following proposition.

\begin{proposition} \label{thm:pjsheaf}
For every object $X$ of $\catE$, $\peejay(X)$ is a sheaf.
\end{proposition} 

\begin{proof} By Proposition~\ref{what-peejay-classifies}, it suffices to
show that every dense monomorphism $m : B \rightarrowtail A$ induces
by pullback an isomorphism 
\[
m^* : \ssubj(X)(A) \rightarrow \ssubj(X)(B)  \, . 
\]
Let us define a proposed inverse $m_\sharp$ as follows. For a family of locally
small closed subobjects $T \rightarrowtail B \times X \rightarrow B$, 
we define 
\[
m_{\sharp}(T) = \cl_{A \times X}(T) \, . 
\]
Here, we view $T$ as a subobject of $A \times X$ via composition 
with the evident 
monomorphism $B \times X \rightarrowtail A \times X$. We need to 
show that the result is a locally small family, and that $m_\sharp$ 
and $m^*$ are mutually inverse.
First, note that $m^*$ is just intersection with $B \times X$. Therefore, 
for a locally small closed family $T \rightarrowtail B \times X \rightarrow 
B$ we have
\begin{eqnarray*}
m^* m_\sharp(T) & = & m^* \cl_{A \times X}(T) \\
 & = & \cl_{B \times X}(m^* T) \\
 & = & \cl_{B \times X}(T) \\
 & = & T \, ,
\end{eqnarray*}
where we used the naturality of the closure operation, that
$T \leq A \times X$, and the assumption that $T$ is closed.
This also shows that $m_{\sharp}(T)$ becomes locally small over $B$ when pulled 
back along the dense monomorphism $m : B \rightarrowtail A$. By 
Proposition~\ref{lemma:J-smallness-is-local-2}, it is locally small 
over $A$, as required. Finally, for a locally small closed family
$S \rightarrowtail A \times X \rightarrow A$, we have
 \begin{eqnarray*}
 m_{\sharp} m^*(S) & = & 
m_{\sharp}( S \cap (B \times X) ) \\
& = & \cl_{A \times X}(S \cap (B \times X)) \\
    & = & \cl_{A \times X}(S) \cap \cl_{A \times X}(B \times X) \\
    & = & S \cap (A \times X) \\
    & = & S \, , 
\end{eqnarray*}
as desired.
\end{proof}

\begin{lemma} \label{thm:clsubsh}
A subobject of a sheaf is a sheaf if and only if it is closed.
\end{lemma}

\begin{proof} Let $i : S \rightarrowtail X$ be a monomorphism, and assume
that $X$ is a sheaf. We begin by proving that if $S$ is closed, then it 
is a sheaf. For this, assume given a dense monomorphism $m : 
B \rightarrowtail A$ and a map $f : B \rightarrow S$. 
Define $g : B \rightarrow X$ to be $i \, f : B \rightarrow X$. By the 
assumption that $X$ is a sheaf, there exists a unique $\bar{g} : A
\rightarrow X$ making the following diagram commute
\[
\xymatrix{
B \ar[r]^f \ar@{ >->}[d]_m & S \ar@{ >->}[r]^i & X  \\
A \ar@/_1pc/@{.>}[urr]_{\bar{g}}     &  & }
\]
We have 
\begin{eqnarray*}
\im(\bar{g}) & = & \bar{g}_{!}(A) \\
             & = & \bar{g}_{!}(\cl_A(B)) \\
             & \leq & \cl_X( \bar{g}_{!}(B)) \\
             & = & \cl_X ( \im(g)) \\
             & \leq & \cl_X(S) \\
             & = & S \, ,
\end{eqnarray*}
where the first inequality follows from adjointness. 
Therefore, $\bar{g} : A \rightarrow X$ factors through $\bar{f} : 
A \rightarrow S$, extending  $f : B \rightarrow S$
as required. Uniqueness of this extension follows by the uniqueness
of $\bar{g} : A \rightarrow X$ and the assumption that $i : S
\rightarrow X$ is a monomorphism. 

For the converse implication, we need to show that if $S$ is a sheaf 
then it is closed. The monomorphism 
$m : S \rightarrowtail \cl_X(S)$ is dense and therefore the identity 
$1_S : S \rightarrow S$ has an extension $n : \cl_X(S) \rightarrow S$
with $n \, m = 1_S$. But by the preceding part, we know that $\cl_X(S)$ 
is a sheaf and that $m \, n \, m = n = 1_{\cl_X(S)} \cdot m$. Therefore,
both $n \cdot m$ and $1_{\cl(S)}$ are extensions of $m : S \rightarrowtail
\cl_X(S)$ to $\cl_X(S)$. Since $S$ is a sheaf, they must be equal. 
Thus $m$ and $n$ are mutually inverse, and so $S \iso \cl_X(S)$ as desired.
\end{proof}

Lemma~\ref{thm:pjsheaf} and Lemma~\ref{thm:clsubsh} imply that
the object $\asf(X)$ is a sheaf  for every $X \in \catE$,
since $\asf(X)$ is a closed subobject of the sheaf $\peejay(X)$. 
In order to show that we have indeed
defined a left adjoint to the inclusion, we need the following
Lemma~\ref{thm:next},
concerning the map $\eta_X : X \rightarrow \asf(X)$, which
was defined at the beginning of the section, using the map $\sigma_X : X \rightarrow 
\peejay(X)$ taking $x : X$ into $[ \{ x \} ] : \peejay(X)$.

\begin{lemma} \label{thm:next} For every $X \in \catE$, we have
\[
\myker(\eta_X) = \myker(\sigma_X) = \cl_{X \times X}(\Delta_X)
\]
in $\sub(X \times X)$.
\end{lemma}

\begin{proof} For $x,y:X$, we have
\begin{eqnarray*}
\sigma(x) = \sigma(y) & \Iff & \cl(\{ x \}) = \cl(\{ y \}) \\
& \Iff & 
(\exists p : \myomega)\big( J(p) \land \big( 
p \Rightarrow x \in \{ y \} \big) \big)  \\
 & \Iff & (x, y) \in \cl_{X \times X}(\Delta_X) \, , 
\end{eqnarray*}
as required.
\end{proof}

Generally, a map $f:X \rightarrow Y$ with $\myker(f) \leq \cl(\Delta_X)$ is 
called \emph{codense}.

\begin{theorem}\label{thm:uniprop} \label{thm:asf}
For every $X \in \catE$, $\asf(X) \in \JSh$ is the associated sheaf
of $X$, in the sense that for every map $f : X \rightarrow Y$ into a
sheaf $Y$, there exists a unique $\bar{f} : \asf(X) \rightarrow Y$ 
making the following diagram commute
\[
\xymatrix{
X \ar[r]^(.43){\eta_X} \ar@/_1pc/[dr]_f & \asf(X) \ar[d]^{\bar{f}} \\
                  & Y }
\]
The resulting left adjoint $\asf : \catE \rightarrow \JSh$ preserves
finite limits.
\end{theorem}

\begin{proof}  Given $f : X \rightarrow Y$,  $\myker(f)$ 
is a pullback of $Y \rightarrowtail Y \times Y$. Since $Y \times Y$ 
is a sheaf, Lemma~\ref{thm:clsubsh} implies that $Y$ is closed, and 
hence $\myker(f)$ is closed in $X \times X$. But 
certainly $\Delta_X \leq \myker(f)$ and therefore
\[
\myker(\sigma_X) = \cl(\Delta_X) \leq \myker(f) \, . 
\] 
Since any epimorphism is the coequaliser of its kernel pair, $f$ factors 
uniquely through the codense epimorphism $X \twoheadrightarrow 
\im(\sigma_X)$. Since $Y$ is a sheaf, the map from $\im(\sigma_X)$ to $Y$ 
extends uniquely along 
the dense monomorphism $\im(\sigma_X) \rightarrowtail \asf(X)$, giving a 
unique factorisation of $f$ through $\eta_X$ as desired.

To show that the associated sheaf functor preserves finite limits, we may 
proceed exactly as 
in~\cite[\S V.3]{MacLaneS:shegl}, since the argument there uses only the 
structure of a Heyting category on $\catE$ and  the fact that the associated 
sheaf functor is defined by embedding each object $X$ by a codense map into 
a sheaf. In particular, injectivity of the sheaf in which $X$ is embedded 
is not required to carry over the proof. 
\end{proof}

Theorem~\ref{thm:uniprop} allows us to deduce 
Proposition~\ref{thm:sheaves-is-hpretopos}, which 
contains the first part of Theorem~\ref{thm:main}.
After stating it, we discuss in some detail the structure of
the category $\JSh$. 

\begin{proposition} \label{thm:sheaves-is-hpretopos}
The category $\JSh$ is a Heyting pretopos.
\end{proposition}

\begin{proof} This is an immediate consequence of Theorem~\ref{thm:uniprop}.
\end{proof}

We discuss the relationship between the structure of $\catE$ and that
of $\JSh$ is some detail. Here and subsequently, operations in $\catE$,
will be denoted without subscript, such as $\im$, while their counterparts 
in $\JSh$ will be denoted with subscript, such as $\im_J$. Limits 
in $\JSh$ are just limits in $\catE$, since sheaves are closed under 
limits. Colimits are the sheafifications of colimits taken in $\catE$, 
since the associated sheaf functor, being a left adjoint, preserves colimits. 
The lattice $\sub[J](X)$ is the  sub-lattice of closed elements of 
$\sub(X)$. Here, the closure operation is a reflection and therefore 
meets in $\sub[J](X)$ are meets in $\sub(X)$, and joins in $\sub[J](X)$
are closures of joins in $\sub(X)$. Moreover, for $S,T \in \sub(X)$, 
we have
\begin{eqnarray*} S \land \cl(S \Implies T) & \leq & 
\cl(S) \land \cl(S \Implies T) \\
    & = & \cl (S \land (S \Implies T)) \\
    & \leq & \cl(T) \\
    & = &  T \, . 
\end{eqnarray*}
Therefore $\cl (S \Implies T) \leq (S \Implies T)$ and so $S \Implies T$ is 
closed. Thus, the implication of $\sub(X)$ restricts so as to give 
implication in $\sub[J](X)$. Since limits in $\catE_J$  agree with limits in $\catE$, 
the inverse image functors in $\JSh$ are just the restrictions of the 
inverse image functors of $\catE$. 

It is clear that the sheafification of a dense monomorphism 
$m: B \rightarrowtail A$ is an isomorphism, since for any sheaf 
$X$ we have
\[
\catE_J\big(\asf(A),X \big) \iso \catE(A, X) \iso \catE(B, X) 
\iso \catE_J \big( \asf(B), X \big) \, .
\]
Therefore, epimorphisms in~$\JSh$ are precisely the dense maps in 
$\catE$, since 
coequalisers in $\JSh$ are sheafifications of coequalisers in $\catE$, 
and dense maps are the sheafifications of 
epimorphisms of $\catE$. We have seen before that dense maps are stable 
under pullback. Consequently, quasi-pullbacks in $\JSh$ are exactly 
those squares that are local quasi-pullbacks in $\catE$.

If $(r_1,r_2): R \rightarrowtail X \times X$ is an equivalence relation in 
$\JSh$, it is also an equivalence relation in $\catE$, so has 
an effective quotient $q: X \twoheadrightarrow Q$ in $\catE$. Its 
sheafification $\asf (q) = \eta_Q \cdot q : X \to \asf(Q)$ is then a 
coequaliser for $(r_1,r_2)$ in sheaves; but this is an effective 
quotient for $R$, since 
\[
\myker(\eta_Q \cdot q) = q^* \cl(\Delta_Q) = \cl (\myker q) = \cl (R) \, ,
\] 
and $R$ is closed in $X \times X$, as it is a subsheaf. 

Images in $\JSh$ are the closures of images in $\catE$. For 
$f: X \rightarrow Y$, $\im_J(f) \rightarrowtail Y$ is the least 
subsheaf of $Y$ throught which $f$ factors, so the least closed 
subobject of $Y$ containing $\im(f)$, which can be readily identified
with $\cl_Y(\im f)$. 
Consequently, the corresponding forward image functor
\[
(\exists_{J})_{f} : \sub[J](X) \rightarrow \sub[J](Y)
\] 
is given by 
\[
(\exists_J)_f(S) = \cl_Y(\exists_f S) \, . 
\] 
Dual image functors in $\JSh$ are again just the restrictions of 
the dual images functors of $\catE$, since if $S \rightarrowtail X$
is closed so is $\forall_f(S) \rightarrowtail Y$, and therefore 
$\cl (\forall_f S) \leq \forall_f (\cl S)$. These give us an immediate 
translation from the internal logic of $\JSh$ into that of $\catE$.

\section{Small maps in sheaves}
\label{sec:smams}

 Within this section, we study 
locally small maps between sheaves, completing the proof of 
Theorem~\ref{thm:main}. The associated sheaf functor followed
by the inclusion of sheaves into $\catE$ preserves dense maps.
Indeed, avoiding explicit mention of the inclusion functor, 
if~$h : B \rightarrow A$ is dense then we have
\[
\im(\eta_A \, h) = \im ( \asf(h) \, \eta_B) \leq \im( \asf(h)) \, .
\]
But $\eta_A \, h : B \rightarrow \asf(A)$ is a dense map, and so its
image is dense in $\asf(A)$. Also, if $X$  is a sheaf then $\eta_X : 
X \rightarrow \asf(X)$ is easily seen to be an isomorphism. 
 
\begin{lemma} \label{thm:natjquasi}
For every $f : X \rightarrow Y$ in $\catE$, the naturality 
square 
\[
\xymatrix{
X \ar[r]^(.45){\eta_X} \ar[d]_f & \asf(X)  \ar[d]^{\asf(f)} \\
Y \ar[r]_(.45){\eta_Y} & \asf(Y) }
\]
is a local quasi-pullback. 
\end{lemma}

\begin{proof}  The canonical 
map $X \rightarrow \asf X \times_{\asf Y} Y$, which we wish to show to
be dense, factors through $X \times_{\asf Y} Y$ as in the diagram below:
\medskip
\[
\xymatrix{
X    \ar[dr] 
     \ar@/^1pc/[drr] 
     \ar@/^2pc/[drrr]^f 
     \ar@/_1pc/[ddr]_{1_X}   &            
                       &   
                       & 
                       \\
                       &  
X \times_{\asf(Y)} Y          \ar[d] \ar[r]    & 
\asf(X) \times_{\asf(Y)} Y    \ar[d] \ar[r]    & 
Y                             \ar[d]^{\eta_Y}  \\
                                               &  
X                             \ar[r]_{\eta_X}  & 
\asf(X)                       \ar[r]_{\asf(f)}  & 
\asf(Y)                                      }
\]
The map $ X \times_{\asf Y} Y \rightarrow \asf X \times_{\asf Y} Y$ is a 
pullback of $\eta_X$ and therefore it is dense. We just need to show that 
$X \rightarrow X \times_{\asf Y} Y$ is dense. For this, let us consider
the following diagram, in which each square is a pullback
\[
\xymatrix{
X \ar[r]^{\iso}\ar@/_1pc/[dr]_{(f,f)} & \Gamma_f \ar[r] \ar[d] 
 & X \times_{\asf(Y)} Y \ar[r] \ar[d] & 
X \ar[d]^f \\
 & \Delta_Y \ar[r] \ar@/_1pc/[dr]_{\iso} & 
\cl(\Delta_Y) \ar[r]^(.53){\pi_1} \ar[d]^{\pi_2} & Y \ar[d]^{\eta_Y} \\
 &                        & Y \ar[r]_{\eta_Y} & \asf(Y) }
\]
We then have the following chain of isomorphisms of subobjects 
of $X \times Y$
\begin{eqnarray*}
  X \times_{\asf Y}Y & \iso & X \times_Y \myker(\eta_Y) \\
                    & \iso & (f,1_Y)^* \cl_{Y \times Y} (\Delta_Y) \\
                    & \iso & \cl_{X \times Y}\big((f,1_Y)^* \Delta_Y \big) \\
                    & \iso & \cl_{X \times Y}(\Gamma_f) \, . 
\end{eqnarray*}
Therefore, the map $X \rightarrow X \times_{\asf Y} Y$ can be identified as
the map $\Gamma_f \rightarrow \cl(\Gamma_f)$, which is dense.
\end{proof}

\begin{corollary} \label{cor:j-smallness-from-asf} \hfill
\begin{enumerate}[(i)]
\item If $f : X \rightarrow Y$ is small, then $\asf(f) : 
\asf(X) \rightarrow \asf(Y)$ is locally small.
\item If $f: X \rightarrow A$ is locally small, then there is some small 
map $g : Y \rightarrow B$ and dense $h: \asf(B) \rightarrow A$ 
such that the following diagram is a pullback
\[
\xymatrix{
\asf(Y) \ar[r] \ar[d]_{\asf(g)} & X \ar[d]^{f} \\
\asf(B) \ar[r]_h                & A }
\]
\end{enumerate}
\end{corollary}

\begin{proof} Since the components of the unit of the adjunction
are dense maps, part $(i)$ follows from Corollary~\ref{thm:addthis} and
Lemma~\ref{thm:natjquasi}. For part $(ii)$, use the definition
of locally small map to get a dense map $h : B \rightarrow A$ and an
indexed family of subobjects $Y \rightarrowtail B \times X \rightarrow B$,
together with a diagram of form
\[
\xymatrix{
Y \ar[r] \ar@{ >->}[d] & X \ar@{ >->}[d] \\
B \times X \ar[r] \ar[d] & A \times X \ar[d] \\
B \ar[r]_h & A }
\]
such that $Y \rightarrowtail B \times_A X$ is dense. Then, let $g : Y \rightarrow B$ be the evident map, which is small. The required square can be readily
obtained by recalling that $\asf : \catE \rightarrow \JSh$ preserves 
pullbacks, sends dense monomorphisms to isomorphisms, and preserves dense 
maps. 
\end{proof}

Corollary~\ref{cor:j-smallness-from-asf} allows us to regard the family of
locally small maps as the smallest family of maps in $\JSh$ containing 
the sheafifications of the small maps in $\catE$ and closed under descent 
along dense maps.

\begin{lemma} \label{lemma:compositions-j-small} \hfill 
\begin{enumerate}[(i)] 
\item Identity maps are locally small.
\item Composites of locally small maps are locally small.
\end{enumerate}
\end{lemma}

\begin{proof} All identities are trivially locally small. For composition, 
suppose $f: X \rightarrow A$ and $f' : X' \rightarrow X$ are locally small. We construct the diagram
\[
\xymatrix{
Y' \times_C Z \ar[r] \ar[d] \ar@{}[dr]|{(5)} & 
Z \ar[r] \ar[d]^h \ar@{}[dr]|{(3)} &  
Y \times_X X' \ar[r] \ar[d]^{f''} \ar@{}[dr]|{(2)} & 
X' \ar[d]^{f'} \\
Y' \ar[r] \ar[d] \ar@{}[drr]|{(4)} & 
C \ar@{->>}[r]        & 
Y \ar[r] \ar[d]^{g} \ar@{}[dr]|{(1)} & 
X \ar[d]^f \\
B' \ar@{->>}[rr]       &     
              & 
B \ar@{->>}[r]_{p}        & 
A}
\]
as follows. We begin by applying Lemma~\ref{j-small-as-j-q-pb} to $f : X 
\rightarrow A$ so as to obtain the local quasi-pullback (1), where
$g : Y \rightarrow B$ is a small map and $p : B \rightarrow A$ is
an epimorphism.
Then, we construct the pullback (2) and obtain the locally small map 
$f'' : Y \times_X X' \rightarrow Y$. We can apply 
Lemma~\ref{j-small-as-j-q-pb} to it so as to obtain the 
local quasi-pullback in~(3), where 
$h : Z \rightarrow C$ is a small map and 
$C \twoheadrightarrow Y$ is an epimorphism. Next, we apply
the Collection Axiom to the small map $Y \rightarrow B$ and
the epimorphism $C \twoheadrightarrow Y$, so as to obtain (4), where
$Y' \rightarrow B'$ is a small map and $B' \twoheadrightarrow B$ is
an epimorphism. Finally, (5) is a pullback and therefore 
$Y' \times_C Z \rightarrow Y'$ is small since $h : Z \rightarrow C$
is so. To conclude that the composite of $f' : X' \rightarrow X$ and
$f : X \rightarrow A$ is locally small, it is sufficient to 
apply Lemma~\ref{j-small-as-j-q-pb} to the whole diagram,
which is a local quasi-pullback since it is obtained as the pasting
of local quasi-pullbacks. 
\end{proof}

\begin{lemma} \label{thm:quotients-in-sheaves} For every commutative 
diagram of the form
\[
\xymatrix{
X \ar[rr]^{d} \ar[dr]_{f} &   & X' \ar[dl]^{f'} \\
                  & A & }
\]
if  $f : X \rightarrow A$ is locally small and $d : X \rightarrow X'$ is 
dense, then $f' : X' \rightarrow A$ is locally small. 
\end{lemma}

\begin{proof} By Lemma~\ref{j-small-as-j-q-pb} we have a 
local quasi-pullback 
\[
\xymatrix{
 Y \ar[r]^{k} \ar[d]_{g} & X \ar[d]^f \\
 B \ar@{->>}[r]_{h} & A }
\]
where $g : Y \rightarrow B$ is a small map and $h : B \twoheadrightarrow A$ is 
an epimorphism. Since $f : X \rightarrow A$ is 
$f'\, d : X \rightarrow A$, we can expand the diagram above as follows
\[
\xymatrix{
Y   \ar@/^1pc/[drr]^{k} \ar[dr] \ar@/_1pc/[ddr] \ar@/_2pc/[dddr]_{g} &   &   \\
    & B \times_A X  \ar[d] \ar[r] & X \ar[d]^{d}   \\     
    & B \times_A X' \ar[d] \ar[r] & X' \ar[d]^{f'} \\
    & B             \ar[r]_{h}        & A }
\]
The map $Y \rightarrow B \times_A X'$ is dense since it is the composite
of the dense map $Y \rightarrow B \times_A X$ with the pullback of the
dense map $d : X \rightarrow X'$. Hence, the commutative square
\[
\xymatrix{
 Y \ar[r]^(.43){d \, k} \ar[d]_{g} & X' \ar[d]^{f'} \\
 B \ar[r]_{h} & A }
\]
is a local quasi-pullback, and so $f' : X' \rightarrow A$ is locally
small. 
\end{proof}

\begin{lemma} \label{lemma:strong-collection-in-sheaves} Let
$X, A, P$ be sheaves.  For every locally
small map $f : X \rightarrow A$ and every dense map $P \rightarrow X$, 
there exists a local quasi-pullback of the form
\[
\xymatrix{
Y \ar[r] \ar[d]_g & P \ar[r] & X \ar[d]^f \\
B \ar[rr]_h       &          & A }
\]
where $g : Y \rightarrow B$ is locally small and $h : B \rightarrow A$ is
dense.
\end{lemma}

\begin{proof}Given such a pair, we construct the diagram
\[
\xymatrix{
           &        &          &   & Y \times_X P   \ar[r]  \ar[dd]_e 
\ar@{}[ddr]|{(2)}   & 
P \ar[dd]^{d} \\
Z' \times_C Z \ar[r] \ar[d]  \ar@{}[dr]|{(6)} & Z \ar[r] \ar[d] 
\ar@{}[drr]|{(4)}
& Y \times_X P \ar[r] \ar@/^1pc/[urr]^{1_{Y \times_X P}} 
 & Y' \ar[d] &             & \\
Z' \ar[r] \ar[d]  \ar@{}[drrrr]|{(5)}     & C \ar[rr]       &       & 
B'  \ar[r] \ar@{}[uur]|{(3)} & 
Y \ar[r] \ar[d] \ar@{}[dr]|{(1)}  & X \ar[d]^{f} \\
C' \ar[rrrr]     &                   &            & & B \ar[r]  & A}
\]
as follows. The commutative square (1) is obtained by applying 
Lemma \ref{j-small-as-j-q-pb} to $f : X \rightarrow A$. In 
particular, it is a local quasi-pullback. Diagram (2) is a pullback.
To construct (3), first we factor the dense map $e : Y \times_X P 
\rightarrow Y$ first
as an epimorphism $Y \times_X P \rightarrow Y'$ followed by a dense 
monomorphism $Y' \rightarrow Y$, and then we factor the dense monomorphism as a small
dense monomorphism $Y' \rightarrow B'$ followed by an epimorphism
$B' \rightarrow Y$ by Lemma~\ref{lemma:handy-factorisation}. We can
then apply the Collection Axiom in $\catE$ to construct a diagram
in (4) which is a quasi-pullback. By the definition of quasi-pullback
and the fact that $Y' \rightarrow B'$ is dense, it follows that 
$Z \rightarrow C$ is a dense map since it is the composition
of a dense map with an epimorphism. We apply again the Collection
Axiom to construct (5), and finally (6) is obtained by another 
pullback.

Since the pasting of (1) and (5) is a local quasi-pullback and the map
$Z' \times_C Z \rightarrow Z'$ is dense, the resulting diagram
\[
\xymatrix{
Z' \times_C Z \ar[r] \ar[d] & P \ar[r]^{d} & X \ar[d]^{f} \\
C' \ar[rr] &  & A }
\]
is a local quasi-pullback as well. Note that some of the objects in
the diagram above need not be sheaves, since they have been 
obtained by applying the Collection Axiom in $\catE$. To complete 
the proof, it suffices to apply the associated sheaf functor, so
as to obtain the diagram
\[
\xymatrix{
\asf(Z' \times_C Z) \ar[r] \ar[d] & P \ar[r]^{d} & X \ar[d]^{f} \\
\asf(C') \ar[rr] &  & A }
\]
This provides the required diagram, since the associated sheaf functor
preserves dense maps and pullbacks and sends small maps into locally
small maps.
\end{proof}

\begin{proposition} \label{thm:small-maps-in-sheaves}
The family $\mapS_J$ of locally small maps in $\JSh$ satisfies
the axioms for a family of small maps.
\end{proposition}

\begin{proof} Lemma~\ref{lemma:compositions-j-small} proves Axiom (A1). 
Axiom~(A2), asserting stability under pullbacks, holds by the very
definition of locally small map, as observed before 
Definition~\ref{thm:jsmall}. Axiom~(A3) follows by 
Proposition~\ref{lemma:J-smallness-is-local}. 
Axioms (A4) and (A5) follow from the corresponding axioms in $\catE$, using
the fact that sheafification of small maps is locally small and the fact that 
the initial object and the coproducts of $\JSh$ are the sheafification of 
the initial object and of coproducts in $\catE$, respectively. 
Lemma~\ref{thm:quotients-in-sheaves} proves Axiom (A6), and 
Lemma~\ref{lemma:strong-collection-in-sheaves} proves Axiom (A7). 
Axiom~(P1) follows by Proposition~\ref{what-peejay-classifies}
and Lemma~\ref{thm:clsubsh}.
\end{proof}

\begin{remark} The theorem asserting that the category internal sheaves over a 
Lawvere-Tierney local operator in an elementary topos is again an
elementary topos~\cite[\S V.2]{MacLaneS:shegl} can be seen as following
from a special case of Proposition~\ref{thm:sheaves-is-hpretopos} and 
Proposition~\ref{thm:small-maps-in-sheaves}. Indeed, elementary toposes 
are examples of our setting, as explained in Example~\ref{exa:topos}, and 
Lawvere-Tierney coverages are equivalent to Lawvere-Tierney local operators, 
as explained in Remark~\ref{rem:lawvere-tierney}. By 
Proposition~\ref{thm:sheaves-is-hpretopos}, we know that the category of 
sheaves is a Heyting pretopos,
and by Proposition~\ref{thm:small-maps-in-sheaves} we know that it has a 
subobject classifier. 
\end{remark}

\begin{proposition} \label{thm:wrep-in-sheaves} 
If the family of small maps $\mapS$ in $\catE$ satisfies the
Exponentiability and Weak Representability Axiom, so does the
family of locally small maps $\mapS_J$ in $\catE_J$.
\end{proposition}

\begin{proof} Since we know by Proposition~\ref{thm:small-maps-in-sheaves}
that axiom (P1) holds for $\mapS_J$ and (P1) implies (S1), 
as recalled in Section~\ref{sec:asthp}, we have that (S1) holds for
$\mapS_J$. For (S2), we need to show that there exists a locally small 
map $u_J : E_J 
\rightarrow U_J$ in $\JSh$ such that every locally small map 
$f : X \rightarrow A$ fits into a diagram of form
\[
\xymatrix{
X \ar[d]_f & Y \ar[d] \ar[r] \ar[l] & E_J \ar[d]^{u_J} \\  
A        & B \ar[r] \ar@{->>}[l]^h        & U_J }
\]
where the left-hand square is a quasi-pullback, the right-hand side
square is a pullback and $h : B \twoheadrightarrow A$ is an epimorphism. 
Let $u : E \rightarrow U$ the weakly representing map
for small maps in $\catE$, which exists by our assumption
that $\mapS$ satisfies (S2). We define 
$u_J : E_J \rightarrow U_J$ 
to be its sheafification, $\asf(u) : \asf(E) \rightarrow \asf(U)$.
Given $f : X \rightarrow
A$ locally small, let us consider a diagram 
\begin{equation}
\label{equ:recall}
\xymatrix{
Y \ar[d]_g \ar[r]^h & X \ar[d]^f \\
B \ar[r]_h & A }
\end{equation}
as in part $(ii)$ of Corollary~\ref{cor:j-smallness-from-asf}. In particular,
$g : Y \rightarrow B$ is small in $\catE$. By the Weak Representability
Axiom for small maps applied to $g : Y \rightarrow B$, we obtain a diagram 
of form
\[
\xymatrix{
X \ar[d]_f & \asf(Y) \ar[l] \ar[d] & \asf(Z) \ar[r] \ar[d] 
\ar[l] & \asf(E) \ar[d]^{\asf(u)}  \\
A        & \asf(B)  \ar[l]       & \asf(C) \ar[r] \ar[l]        & \asf(U) }
\]
This is obtained by applying the associated sheaf functor to the diagram
expressing that $g : Y \rightarrow B$ is weakly classified
by $u : E \rightarrow U$ and to the one in~(\ref{equ:recall}), recalling
that $X \iso \asf(X)$ if $X$ is a sheaf. The desired conclusion follows 
by recalling that the associated sheaf functor preserves dense maps and 
pullbacks. 
\end{proof}

The combination of Proposition~\ref{thm:sheaves-is-hpretopos},
Proposition~\ref{thm:small-maps-in-sheaves}, and 
Proposition~\ref{thm:wrep-in-sheaves}
provides a proof of Theorem~\ref{thm:main}.

\begin{remark} We may note that for a locally small sheaf $X$, witnessed by
an epimorphism $h : B \twoheadrightarrow 1$ and an indexed family of small
subobjects $S \rightarrowtail B \times X \rightarrow B$, fitting in a diagram
\[
\xymatrix{
S \ar[r]^{k} \ar@{ >->}[d] & X \ar@{ >->}[d] \\
B \times X \ar[r] \ar[d] & 1 \times X \ar[d] \\
B \ar@{->>}[r]_{h} & 1 }
\]
with $S \rightarrowtail B \times X$ dense, an exponential of the form 
$Y^X$, where $Y$ is a sheaf, can be constructed explicitly as the 
quotient of 
\[
\sum_{b:B} Y^{S_b}
\] 
by the equivalence relation defined by
\[ 
(b,f) \sim (b',f')
\Iff 
(\forall s  : S_{b}) 
(\forall s' : S_{b'}) 
\left[ k(s) = k(s')  \Implies\ f(s) = f'(s') \right] \, .
\]
\end{remark}

\begin{remark}
Of course, if $\catE$ has more structure than we have assumed so far, or 
satisfies stronger axioms, we would like these to be preserved by the 
construction of internal sheaves. For example, it is immediate to see
that if~$\catE$ has small diagonals, so does $\JSh$, and that existence
of a universal object in~$\catE$ in the sense of~\cite{SimpsonA:eleacc}
implies the existence of a universal object in $\JSh$. 
\end{remark}

\section*{Acknowledgements}

We wish to thank the organizers of the Summer School on Topos Theory in 
Haute-Bodeux (Belgium), where the idea for the research described here was 
first discussed.  Some of the research was conducted during visits of the
second author to Carnegie Mellon University, for which he wishes to 
express his gratitude. He also wishes to thank the Centre de Recerca
Matem\`atica (Barcelona) for its support.


\providecommand{\bysame}{\leavevmode\hbox to3em{\hrulefill}\thinspace}
\providecommand{\MR}{\relax\ifhmode\unskip\space\fi MR }
\providecommand{\MRhref}[2]{%
  \href{http://www.ams.org/mathscinet-getitem?mr=#1}{#2}
}
\providecommand{\href}[2]{#2}

\end{document}